\documentclass[12pt]{amsart}
\usepackage{amssymb}
\usepackage[all]{xy}

\long\def\forget#1\forgotten{}
\newcommand{\issuenumber}{17}
\newcommand{\issuemonth}{June}
\newcommand{\issueyear}{2006}

\newtheorem{thm}{Theorem}[section]
\newtheorem{prob}[thm]{Problem}

\newtheorem{issue}{Issue}

\theoremstyle{definition}

\theoremstyle{remark}

\newcommand{\ed}{
\general\end{document}}

\newcommand{\Cantor}{{\{0,1\}^\N}}
\newcommand{\oo}{\infty}

\newcommand{\fb}{\mathfrak{b}}

\newcommand{\fd}{\mathfrak{d}}
\newcommand{\fp}{\mathfrak{p}}

\newcommand{\NON}{{\mathsf   {NON}}}
\newcommand{\COF}{{\mathsf   {COF}}}

\newcommand{\M}{\mathcal{M}}

\newcommand{\cov}{\mathsf{cov}}

\newcommand{\cf}{\mathsf{cf}}

\newcommand{\R}{\mathbb{R}}

\newcommand{\fo}{\mathfrak{od}}

\newcommand{\ft}{\mathfrak{t}}

\renewcommand{\split}{\mathsf{Split}}
\newcommand{\bq}{\begin{quote}}
\newcommand{\eq}{\end{quote}}
\newcommand{\cO}{\mathcal{O}}
\newcommand{\B}{\mathcal{B}}
\newcommand{\BG}{\B_\Gamma}

\newcommand{\BO}{\B_\Omega}

\newcommand{\sone}{\mathsf{S}_1}    \newcommand{\sfin}{\mathsf{S}_{fin}}

\newcommand{\ufin}{\mathsf{U}_{fin}}

\newcommand{\nin}{\not\in}

\newcommand{\N}{\mathbb{N}}

\newcommand{\sbst}{\subseteq}
\newcommand{\by}[2]{\par\hfill\emph{#1}, #2}
\newcommand{\nby}[1]{\par\hfill\emph{#1}}
\newcommand{\Tau}{\mathrm{T}}
\newcommand{\CE}{\textsc{CE}}

\newcommand{\be}{\begin{enumerate}}
\newcommand{\ee}{\end{enumerate}}
\newcommand{\bi}{\begin{itemize}}
\newcommand{\ei}{\end{itemize}}

\newcommand{\general}{\small\vfill\par\noindent\hrulefill\par
\noindent\textbf{Previous issues.} The first issues of this
bulletin, which contain general information (first issue), basic
definitions, research announcements, and open problems (all
issues) are available online, on \arx{math.GN/$x$}, where $x$ is
\texttt{0301011}, \texttt{0302062}, \texttt{0303057},
\texttt{0304087}, \texttt{0305367}, \texttt{0312140},
\texttt{0401155}, \texttt{0403369}, \texttt{0406411},
\texttt{0409072}, \texttt{0412305}, \texttt{0503631},
\texttt{0508563}, \texttt{0509432}, \texttt{0512275}, and
\texttt{0603290},
respectively, for issues number $1$ to $16$.\\[0.1cm]
\textbf{Contributions.}
Please submit your contributions (announcements, discussions, and open problems)
by e-mailing us. It is preferred to write them
in \LaTeX{}.
The authors are urged to use as standard notation as possible, or otherwise give
the definitions or a reference to where the notation is explained.
Contributions to this bulletin would not require any transfer of copyright,
and material presented here can be published elsewhere.\\[0.1cm]
\textbf{Subscription.}
To receive this bulletin (free) to your
e-mailbox, e-mail us:\\
{boaz.tsaban@weizmann.ac.il}
}

\newcommand{\nArxPaper}[5]{\subsection{#2}{#4}\par\hfill{\arx{#1}}\par\hfill\emph{#3}}
\newcommand{\DOIpaper}[5]{\subsection{#2}{#4}\par\hfill{\texttt{http://dx.doi.org/#1}}\par\hfill\emph{#3}}

\newcommand{\nAMSPaper}[5]{\subsection{#2}{#4}\par\hfill{\texttt{#1}}\par\hfill\emph{#3}}

\newcommand{\arx}[1]{\texttt{http://arxiv.org/abs/#1}}
\newcommand{\url}[1]{\bq\texttt{#1}\eq}
\newcommand{\online}[1]{The paper is available online at \url{#1}}

\title[$\mathcal{SPM}$ Bulletin \textbf{\issuenumber} (\issuemonth{} \issueyear)]{%
$\mathcal{SPM}$ Bulletin\\[0.5cm]
Issue number \issuenumber: \issuemonth{} \issueyear{} \CE{}}

\begin{document}
\maketitle

\centerline{\emph{This issue is dedicated to Gary Gruenhage}}

\tableofcontents

\section{Editor's note}

Those of you who missed the second SPM meeting in Lecce,
may wish to have a look at Section \ref{SPM}.
In fact, also those who attended may wish to do so.

It is encouraging that each issue comes with so many
interesting announcements of new results in the field
of SPM and in closely related field. This flourishing a good sign.

\medskip

Contributions to the next issue are, as always, welcome.

\medskip

\by{Boaz Tsaban}{boaz.tsaban@weizmann.ac.il}

\hfill \texttt{http://www.cs.biu.ac.il/\~{}tsaban}

\section{Research announcements}

\subsection{Lecce Workshop presentations available online}\label{SPM}
The slides of the talks given at the \emph{Second Workshop on Coverings, Selections, and Games in Topology}
(Lecce, December 2005) are now available online:
\url{http://www.cs.biu.ac.il/\~{}tsaban/SPMC05/slides.html}

\nAMSPaper{http://www.ams.org/journal-getitem?pii=S0002-9939-06-08207-4}
{Borel cardinalities below $c_0$}
{Michael Ray Oliver}
{The Borel cardinality of the quotient of the power set of the natural numbers by the ideal
$\mathcal{Z}_0$ of asymptotically zero-density sets is shown to be the same as that of the equivalence
relation induced by the classical Banach space $c_0$.
We also show that a large collection of ideals introduced by Louveau and Velickovic, with pairwise
incomparable Borel cardinality, are all Borel reducible to $c_0$.
This refutes a conjecture of Hjorth and has facilitated further work by Farah.
}

\nArxPaper{math.GR/0603513}
{Hereditarily non-topologizable groups}
{G\'abor Luk\'acs}
{A group $G$ is non-topologizable if the only Hausdorff group topology that $G$ admits is the discrete one.
Is there an infinite group $G$ such that $H/N$ is non-topologizable for every subgroup $H \le G$
and every normal subgroup $N \lhd H$?
We show that a solution of this essentially group theoretic question provides a solution to the problem of
$c$-compactness.
}

\nArxPaper{math.LO/0603691}
{A hodgepodge of sets of reals}
{Arnold W.\ Miller}
{We prove a variety of results concerning singular sets of reals.
Our results concern: Kysiak and Laver-null sets,
Kocinac and $\gamma_k$-sets,
Fleissner and square $Q$-sets,
Alikhani-Koopaei and minimal $Q$-like-sets,
Rubin and $\sigma$-sets, and Zapletal and the Souslin number.
In particular we show that $\sigma$-sets are Laver-null,
the union of $\gamma_k$-sets need not be $\gamma_k$,
the existence of $Q$-set implies an $\aleph_1$-universal $G_\delta$,
minimal $Q$-like sets which are not $Q$-sets exist,
thin sets need not exist,
and $\mathfrak{sn}^*$ is bounded by the cardinality of the smallest nonmeager set.
}

\nArxPaper{math.LO/0604085}
{Random gaps}
{James Hirschorn}
{It is proved that there exists an $(\aleph_1,\aleph_1)$ Souslin gap in the Boolean algebra
$(L(\nu)/Fin,\subseteq^*_{ae})$ for every nonseparable measure $\nu$.
Thus a Souslin, also known as destructible, $(\aleph_1,\aleph_1)$ gap in $P(\N)/Fin$
can always be constructed from uncountably many random reals.
We explain how to obtain the corresponding conclusion from the hypothesis that Lebesgue measure
can be extended to all subsets of the real line (RVM), and why this runs counter to authoritative
opinions on the nature of consequences of RVM.}

\nArxPaper{math.LO/0604156}
{Covering a bounded set of functions by an increasing chain of slaloms}
{Masaru Kada}
{A slalom is a sequence of finite sets of length omega.
Slaloms are ordered by coordinatewise inclusion with finitely many exceptions.
Improving earlier results of Mildenberger, Shelah and Tsaban, we prove consistency
results concerning existence and non-existence of an increasing sequence of a certain
type of slaloms which covers a bounded set of functions in the Baire space.}

\nArxPaper{math.GN/0604238}
{Baire-one mappings contained in a usco map}
{Ond\v{r}ej F.K.\ Kalenda}
{We investigate Baire-one functions whose graph is contained in a graph of
usco mapping.
We prove in particular that such a function defined on a metric space with values in
$\mathbb{R}^d$ is the pointwise limit of a sequence of continuous functions with graphs
contained in the graph of a common usco map.}

\DOIpaper{10.1016/j.topol.2005.07.015}
{Applications of $k$-covers II}
{A.\ Caserta, G.\ Di Maio, Lj.\ D.R.\ Ko\v{c}inac, and E.\ Meccariello}
{We continue the study of applications of k-covers to some topological constructions,
mostly to function spaces and hyperspaces.}

\nArxPaper{math.GN/0604451}
{Additivity numbers of covering properties}
{Boaz Tsaban}
{This is an invited chapter in the book ``Selection Principles in
Topology'' (Ljubisa Kocinac, ed.), to appear in the book series
\textbf{Quaderni di Matematica}.

The \emph{additivity number} of a topological property (relative to a given space)
is the minimal number of subspaces with this property whose union does not have
the property. The most well-known case is where this number is greater than
$\aleph_0$, i.e. the property is $\sigma$-additive. We give a rather complete survey
of the known results about the additivity numbers of a variety of topological
covering properties, including those appearing in the Scheepers diagram (which
contains, among others, the classical properties of Menger, Hurewicz,
Rothberger, and Gerlits-Nagy). Some of the results proved here were not
published beforehand, and many open problems are posed.}

\nArxPaper{math.LO/0604536}
{Combinatorial images of sets of reals and semifilter trichotomy}
{Boaz Tsaban and Lyubomyr Zdomskyy}
{Using a dictionary translating a variety of classical and modern covering properties into
combinatorial properties of continuous images, we get a simple way to understand the
interrelations between these properties in ZFC and in the realm of the trichotomy axiom
for upward closed families of sets of natural numbers. While it is now known that the answer
to the Hurewicz 1927 problem is positive, it is shown here that semifilter trichotomy implies
a negative answer to a slightly weaker form of this problem.}

\subsection{Another algebraic equivalent of the Continuum Hypothesis}
{Enrico Zoli}
{A renowned theorem due to Erd\H{o}s and Kakutani states that the
Continuum Hypothesis holds if and only if the set of all nonnull
reals is a union of countably  many  Hamel bases. Adapting Erd\H
os and Kakutani's argument, here I show that the Continuum
Hypothesis holds if and only if the set of all transcendental
reals is a union of countably many transcendence bases.}

\nArxPaper{math.LO/0605022}
{A connection between decomposable ultrafilters and possible cofinalities II}
{Paolo Lipparini}
{We use Shelah's theory of possible cofinalities in order to solve a problem about ultrafilters.\\
THEOREM. Suppose that $\lambda$ is a singular cardinal, $\lambda ' < \lambda$,
and the ultrafilter $D$ is $ \kappa $-decomposable for all regular cardinals $ \kappa $ with $
\lambda ' < \kappa < \lambda $. Then $D$ is either $ \lambda $-decomposable, or $ \lambda ^+$-decomposable.\\
We give applications to topological spaces and to abstract logics.}

\nArxPaper{math.GN/0605469}
{Game Approach to Universally Kuratowski-Ulam Spaces}
{A.\ Kucharski and Sz.\ Plewik}
{We consider a version of the open-open game, indicating its connections with universally Kuratowski-Ulam spaces.
We show that: Every I-favorable space is universally Kuratowski-Ulam; If a compact space $Y$ is I-favorable,
then the hyperspace $\exp(Y)$ with the Vietoris topology is I-favorable, and hence universally Kuratowski-Ulam.
Notions of uK-U and uK-U$^*$ spaces are compared.}

\nAMSPaper{http://www.ams.org/journal-getitem?pii=S0002-9939-06-08401-2}
{On the density of Banach spaces $C(K)$ with the Grothendieck property}
{Christina Brech}
{Using the method of forcing we prove that consistently there is a Banach space of continuous
functions on a compact Hausdorff space $C(K)$ with the Grothendieck property and with density
less than the continuum. It follows that the classical result stating that ``no nontrivial complemented
subspace of a Grothendieck space is separable'' cannot be strengthened by replacing ``is separable''
by ``has density less than that of $\ell_\oo$'', without using an additional set-theoretic assumption.
Such a strengthening was proved by Haydon, Levy and Odell, assuming Martin's axiom and the negation of the
continuum hypothesis. Moreover, our example shows that certain separation properties of Boolean algebras
are quite far from the Grothendieck property.
}

\nArxPaper{math.LO/0606021}
{Antichains in partially ordered sets of singular cofinality}
{Assaf Rinot}
{In their paper from 1981, Milner and Sauer conjectured that for any poset $\langle P,\le\rangle$,
if $\cf(P,\le)=\lambda>\cf(\lambda)=\kappa$,
then $P$ must contain an antichain of size $\kappa$.

We prove that for $\lambda>\cf(\lambda)=\kappa$, if there exists a cardinal $\mu<\lambda$ such that $\cov(\lambda,\mu,\kappa,2)=\lambda$,
then any poset of cofinality $\lambda$ contains $\lambda^\kappa$ antichains of size $\kappa$.

The hypothesis of our theorem is very weak and is a consequence of many well-known axioms such as GCH, SSH and PFA.
The consistency of the negation of this hypothesis is unknown.
}

\nArxPaper{math.LO/0606028}
{Bolzano-Weierstrass principle of choice extended towards ordinals}
{W.\ Kulpa, Sz.\ Plewik and M.\ Turza\'nski}
{The Bolzano-Weierstrass principle of choice is the oldest method of the set theory,
traditionally used in mathematical analysis. We are extending it towards transfinite
sequences of steps indexed by ordinals. We are introducing the notions: hiker's tracks,
hiker's maps and statements $P_n(X, Y, m)$; which are used similarly in finite, countable
and uncountable cases. New proofs of Ramsey's theorem and Erd\"{o}s-Rado theorem are presented
as some applications.}

\nArxPaper{math.GN/0606146}
{Nonequality of Dimensions for Metric Groups}
{Ol'ga V.\ Sipacheva}
{An embeddability criterion for zero-dimensional metrizable topological spaces in zero-dimensional metrizable
topological groups is given. A space which can be embedded as a closed subspace in a zero-dimensional metrizable
group but is not strongly zero-dimensional is constructed; thereby, an example of a metrizable group with
noncoinciding dimensions ind and dim is obtained. It is proved that one of Kulesza's zero-dimensional metrizable
spaces cannot be embedded in a metrizable zero-dimensional group.}

\nArxPaper{math.OA/0606168}
{Not all pure states on $B(H)$ are diagonalizable}
{Charles Akemann and Nik Weaver}
{Assuming the continuum hypothesis, we prove that $B(H)$ has a pure state whose restriction to any masa is not pure.
This resolves negatively an old conjecture of Anderson.}

\nArxPaper{math.LO/0404220}
{A comment on $\fp<\ft$}
{Saharon Shelah}
{Dealing with the cardinal invariants ${\mathfrak p}$ and ${\mathfrak t}$ of
the continuum we prove that
\[{\mathfrak m}\geq {\mathfrak p} = \aleph_2\ \Rightarrow\ {\mathfrak t} =
\aleph_2.\]
In other words, if ${\bf MA}_{\aleph_1}$ (or a weak version of this) holds,
then (of course $\aleph_2\le {\mathfrak p}\le {\mathfrak t}$ and)
${\mathfrak p}=\aleph_2\ \Rightarrow\ {\mathfrak p}={\mathfrak t}$.  The
proof is based on a criterion for ${\mathfrak p}<{\mathfrak t}$.}

\section{Problem of the Issue}

A set of reals is \emph{totally imperfect} if it has no perfect subsets.
There are in the literature several examples
of totally imperfect sets of reals satisfying $\ufin(\cO,\Gamma)$, see
\cite{coc2, BaSh01, SFH}.
However, all of them actually satisfy $\sone(\Gamma,\Gamma)$ \cite{wqn}
or at least $\sone(\Gamma,\cO)$ \cite{ideals, SFH}.
Consequently, all of their continuous images satisfy $\sone(\Gamma,\cO)$
and are therefore totally imperfect \cite{coc2}.

\begin{prob}
Could there be a totally imperfect set of reals $X$ satisfying $\ufin(\cO,\Gamma)$
that can be mapped continuously onto $\Cantor$?
\end{prob}

\nby{Tomasz Weiss}

\newpage

\section{Problems from earlier issues}

\begin{issue}
Is $\binom{\Omega}{\Gamma}=\binom{\Omega}{\Tau}$?
\end{issue}

\begin{issue}
Is $\ufin(\Gamma,\Omega)=\sfin(\Gamma,\Omega)$?
And if not, does $\ufin(\Gamma,\Gamma)$ imply
$\sfin(\Gamma,\Omega)$?
\end{issue}

\stepcounter{issue}

\begin{issue}
Does $\sone(\Omega,\Tau)$ imply $\ufin(\Gamma,\Gamma)$?
\end{issue}

\begin{issue}
Is $\fp=\fp^*$? (See the definition of $\fp^*$ in that issue.)
\end{issue}

\begin{issue}
Does there exist (in ZFC) an uncountable set satisfying $\sone(\BG,\B)$?
\end{issue}

\stepcounter{issue}

\begin{issue}
Does $X \nin \NON(\M)$ and $Y\nin\mathsf{D}$ imply that
$X\cup Y\nin \COF(\M)$?
\end{issue}

\begin{issue}[CH]
Is $\split(\Lambda,\Lambda)$ preserved under finite unions?
\end{issue}

\begin{issue}
Is $\cov(\M)=\fo$? (See the definition of $\fo$ in that issue.)
\end{issue}

\begin{issue}
Does $\sone(\Gamma,\Gamma)$ always contain an element of cardinality $\fb$?
\end{issue}

\begin{issue}
Could there be a Baire metric space $M$ of weight $\aleph_1$ and a partition
$\mathcal{U}$ of $M$ into $\aleph_1$ meager sets where for each ${\mathcal U}'\subset\mathcal U$,
$\bigcup {\mathcal U}'$ has the Baire property in $M$?
\end{issue}

\stepcounter{issue} 

\begin{issue}
Does there exist (in ZFC) a set of reals $X$ of cardinality $\fd$ such that all
finite powers of $X$ have Menger's property $\sfin(\cO,\cO)$?
\end{issue}

\begin{issue}
Can a Borel non-$\sigma$-compact group be generated by a Hurewicz subspace?
\end{issue}

\begin{issue}[MA]
Is there an uncountable $X\sbst\R$ satisfying $\sone(\BO,\BG)$?
\end{issue}

\begin{issue}[CH]
Is there a totally imperfect $X$ satisfying $\ufin(\cO,\Gamma)$
that can be mapped continuously onto $\Cantor$?
\end{issue}

\ed